\documentclass[12pt]{amsart}
\usepackage{graphicx}
\usepackage{amssymb, amsmath, amsthm}
 
\newtheorem{lemma}{Lemma}
\newtheorem{cor}[lemma]{Corollary}
 \newtheorem{thm}[lemma]{Theorem}
\newtheorem{defi}[lemma]{Definition}

\newtheorem{propn}[lemma]{Proposition}
\newcommand{\R}{\mathbb R}

\newcommand{\Z}{\mathbb Z}
\begin{document}
%%%%%%%%%%%%%%%%%%%%%%%%%%%%%%%%%%%%%%%%%%%%%%%%%%%%
\title{A characterization of hyperbolic spaces}
%%%%%%%%%%%%%%%%%%%%%%%%%%%%%%%%%%%%%%%%%%%%%%%%%%%
\author{Indira Chatterji and Graham Niblo}
\date{\today}
%%%%%%%%%%%%%%%%%%%%%%%%%%%%%%%%%%%%%%%%%%%%%%%%%%%%%%%%%%%%%%%%%%%%%%%%%%%%%
\begin{abstract}We show that a geodesic metric space is hyperbolic in the sense of Gromov if and only if intersections of balls have bounded eccentricity. In particular, $\R$-trees are characterized among geodesic metric spaces by the property that the intersection of any two balls is always a ball. Both Gromov hyperbolicity and CAT($\kappa$) geometry can be characterised in terms of the geometry of the intersection of balls.\end{abstract}
%%%%%%%%%%%%%%%%%%%%%%%%%%%%%%%%%%%%%%%%%%%%%%%%%%%%%%%%%%%%%%%%%%%%%%%%%%%%%%%%%%%
\maketitle
%\tableofcontents
\section*{Introduction}
It is well known that in an $\mathbb R$-tree the intersection of any two metric balls is itself a metric ball. In this paper we  will show that this is actually a characterization of  $\R$-trees,and that, more generally, the geometry of the intersection of balls encodes information about the curvature of a geodesic metric space. Recall from \cite{Bible}, \cite{CDP}, \cite{GH} or \cite{Gr} that a geodesic metric space is \emph{hyperbolic} (in the sense of Gromov) if there is a constant $\delta\geq 0$ such that for any geodesic triangle, any one side is contained in the $\delta$-neighborhood of the union of the two other sides. We prove the following characterization of hyperbolicity.
\begin{thm}\label{qb}A geodesic metric space $(X,d)$ is hyperbolic  if and only if the intersection of any two metric  balls is at uniformly bounded Hausdorff distance from a ball. \end{thm}
Studying curvature in terms of the geometry of the intersection of metric balls turns out to be very natural and  both Gromov hyperbolicity and the notion of CAT($\kappa$) geometry may be characterised in theses terms (see Section \ref{ext}). 

Tracking constants in the proof of Theorem \ref{qb}, one can show that the hyperbolicity constant depends only on the eccentricity constant. As Pierre Pansu pointed out to us it is a then an easy observation to deduce that the hyperbolicity bound varies linearly with the eccentricity bound so we obtain:
\begin{cor}\label{RtreesChar}
A geodesic metric is an $\R$-tree if and only if the intersection of any two balls is a ball.
\end{cor}
This characterisation of $\mathbb R$-trees was first conjectured in an early version of this paper and Wenger has independently established the conjecture using very different methods, see \cite{W}. Our approach is entirely self-contained. 
The following notion is crucial for our purposes.

\begin{defi}\label{eccentricity}We say that a set $S$ has \emph{eccentricity less than} $\delta$ (for some $\delta\geq 0$) if there is $R\geq 0$ such that 
$$B(c,R)\subseteq S \subseteq B(c',R+\delta)$$
for some $c,c'\in X$. By convention the eccentricity of the empty set is 0. \end{defi}
We shall see that the intersection of balls having uniformly bounded eccentricity is also equivalent to hyperbolicity (Proposition \ref{hyp=>qb} and Lemma \ref{qb=>hyp}). 

The paper is organized as follows. In Section \ref{bigons}  we discuss the geometry of $(1,q)$-quasigeodesics following an idea of Papasoglu in \cite{Papa} and Pomroy  in \cite{Pomroy}, which is an important step in the proof. Section \ref{div} discusses divergence functions and a quantitative version of a theorem in \cite{Papa} and a classical argument implying hyperbolicity. Section \ref{proofs} collects the proofs of Theorem \ref{qb} and Corollary \ref{RtreesChar}. The idea is to show that hyperbolicity is equivalent to intersections of balls having uniformly bounded eccentricity. One technical difficulty lies in the fact that the centre and radius of a ball are not, in general, well defined. When they are an elementary proof can be given as in Section \ref{ext}. Pomroy's work appeared in his Warwick University Masters dissertation but has never been published. In the Appendix we take this opportunity to place his main theorem on the record with our own variation on the proof. The figures may be found at the back of the paper, inserted prior to the bibliography.

{\it Acknowledgements:} The authors are extremely indebted to Pierre Pansu for pointing out to us that our bound on the hyperbolicity constant implied a linear bound via scaling. We would also like to thank Theo Buehler and Deborah Ruoss for pointing out an error in an early draft of this paper and Chris Hruska for the reference \cite{MM}. Finally we would like to thank Hamish Short for his comments on the classical argument  proving that non-linear divergence of geodesics implies hyperbolicity.
%%%%%%%%%%%%%%%%%%%%%%%%%%%%%%%%%%%%%%%%%%%%%%%%%%%%%%%%%%%%%%%%%%%%%%%%%%%%%%%%%%%%%%%%%%%%%%%%%
\section{Bigons in geodesic metric spaces}\label{bigons}In this section we establish preliminary results concerning the geometry of geodesics and quasi-geodesics in a geodesic metric space. We start with the simple observation that  if two geodesics are synchronously far apart, then they are asynchronously at least half as far apart as well.
\begin{lemma}\label{synchronous} Let $\gamma$ and $\gamma'$ be
geodesics with $\gamma(0)=\gamma'(0)=e$. If there exists $t\geq 0$
such that $d(\gamma(t),\gamma'(t))\geq K$  then $d(\gamma(t),
\gamma'(s))\geq K/2$ for all $s$.
\end{lemma}
\begin{proof} 
Suppose that there is an $s\leq t$ such that $d_0=d(\gamma'(s), \gamma(t))<K/2$. Then since $\gamma$ is a geodesic we have $t=d(\gamma(0), \gamma(t))\leq d(e, \gamma'(s))+d(\gamma'(s), \gamma(t))=s+d_0$. Hence $d_0\geq t-s$. But by hypothesis and the triangle inequality
\begin{eqnarray*}
K\leq d(\gamma'(t), \gamma(t))&\leq& d(\gamma'(t), \gamma'(s))+d(\gamma'(s), \gamma(t))\\
&=& t-s+d_0\leq 2d_0 <K.
\end{eqnarray*}
This is a contradiction. Similarly, if there is an $s> t$ such that $d_0=d(\gamma'(s), \gamma(t))<K/2$. Then since $\gamma$ is a geodesic we have $s=d(\gamma'(0), \gamma'(s))\leq d(e, \gamma(t))+d(\gamma(t), \gamma'(s))=t+d_0$. Hence $d_0\geq s-t$. But by hypothesis and the triangle inequality
\begin{eqnarray*}
K\leq d(\gamma'(t), \gamma(t))&\leq& d(\gamma'(t), \gamma'(s))+d(\gamma'(s), \gamma(t))\\
&=& s-t+d_0\leq 2d_0 <K.
\end{eqnarray*}
This is again a contradiction.
\end{proof}
\begin{defi}For a constant $q\geq 0$ a $(1,q)$-quasigeodesic is a continuous map $\gamma:[0,d]\to X$ such  that $\gamma(0)=\gamma'(0)$, $\gamma(d)=\gamma'(d)$ and for all $t\in [0,d]$
$$|t-t'|-q\leq d(\gamma(t),\gamma(t'))\leq |t-t'|+q.$$ 
The points $\gamma(0)$ and $\gamma(d)$ are said to be the endpoints of $\gamma$. A $(1,q)$ quasi-geodeSIc bigon is a pair $\gamma, \gamma'$ of $(1,q)$ quasi-geodesic bigons which have the same end points.
The images of the quasigeodesics $\gamma$ and $\gamma'$ in $X$ are called the \emph{sides} of the $q$-bigon. For $K\geq 0$, we say that a $q$-bigon is $K$-\emph{fat} if there are $s,t$ such that $d(\gamma(s), \gamma'(t))\geq K$ and that it is $K$-\emph{thin} if it is not $K$-fat.\end{defi}

The following remark provides an easy mechanism for constructing $(1,q)$ quasi-geodesics.

\begin{lemma}Let $K\geq 0$ and $x,y$ in $X$. A \emph{$K$-path} from $x$ to $y$ is a continuous path $\mu$ from $x$ to $y$ such that, for any $z=\mu(t)$ for some $t$
$$d(x,z)+d(z,y)\leq d(x,y)+K.$$
Given any point on $p$ on a $(1,K)$-quasi-geodesic with end points $x,y$ we obtain a  $K$-path by taking a broken geodesic $xpy$. 
\end{lemma}

The main result in this section is that if two geodesics start at the same point and travel ``almost parallel'' for long enough, then this ensures the existence  of  $\epsilon$-bigons for any $\epsilon\geq 0$ whose fatness depends on the distance between the two geodesics.  The proof of this follows the outline in \cite{Papa, Pomroy}. More precisely we prove the following.
\begin{thm}\label{meat}Let $(X,d)$ be a geodesic metric space, and choose integers $K_0, K_1, q$ with $1\leq K_0<K_1$, $q\geq 3$. Assume that there are two geodesics $\gamma$, $\gamma'$ and a real number $R>0$ such that $\gamma(0)=\gamma'(0)$ and $d(\gamma(R+r),\gamma'(R+r))\in[K_0,K_1]$ for all $r\in(0,r_0)$. If  $X$ does not contain a $K_0/2$-fat $1/q$-bigon then 
$$r_0\leq (q(K_1-K_0)+1)(2^{qK_1+1}-1)qK_1+1.$$
In words, there is an upper bound on the length for which two geodesics can travel at a controlled distance.
\end{thm}
\begin{proof}For the purpose of the proof we introduce a few definitions. Take $J\in[K_0,K_1]$, we call an element $t\in (0,r_0)$ a $J$-\emph{point} if $d(\gamma(R+t),\gamma'(R+t))\in[J,J+1/q)$. We say that $t$ is an integral $J$-point if $t$ is a positive integer and a $J$-point. We define the $J$-\emph{distance} between two integral $J$-points $t\not=t'$ by
$$d_J(t,t')=\sharp\{\hbox{integral }J-\hbox{points between }t\hbox{ and }t'\}+1$$
and set  $d_J(t,t')=0$ if and only if $t=t'$. Note that since we are assuming that $t, t'$ are integers $d_J(t, t')\leq t'-t$. The interval $[K_0,K_1]$ is covered by $q(K_1-K_0)+1$ disjoint half open  intervals of the form $[J,J+1/q)$, where $J\in (1/q)\Z\cap[K_0,K_1]$.

Given that $X$ does not contain a $K_0/2$-fat $1/q$-bigon we claim the following:

{\it Claim:}  For any two integral $J$-points $t$ and $t'$, if $d_J(t,t')=(2^n-1)qJ$, then
$$d(\gamma(R+t),\gamma'(R+t))<t'-t+J-n/q+1/q.$$
We prove the above claim by induction over $n$: first, for $n=0$ we have $t=t'$ and $d(\gamma(R+t),\gamma'(R+t'))<J+1/q$ because $t=t'$ is a $J$-point. 

So we assume the claim is true for $n=m$ and show it for $n=m+1$. 

Given integral $J$-points $t, t'$ with $d_J(t, t')=(2^{m+1}-1)qJ$ we let $t_1, t_2$ be the integral $J$ points such that 
\[d_J(t,t_1)=d_J(t_2,t')=(2^m-1)qJ.\]
Note that $d_J(t_1,t_2)=qJ$ because $(2^{m+1}-1)qJ=(2^m-1)qJ+qJ+(2^m-1)qJ$, and hence both $d(\gamma(R+t_1),\gamma(R+t_2))$ and $d(\gamma'(R+t_1),\gamma'(R+t_2))$ are greater than or equal to $qJ$. 

Set $P=\gamma(R+t), Q=\gamma'(R+t), P_1=\gamma(R+t_1), Q_1=\gamma'(R+t_1), P_2=\gamma(R+t_2), Q_2=\gamma'(R+t_2), P'=\gamma(R+t'), Q'=\gamma'(R+t')$, as shown in Figure 1.

%
%\vspace{.2in}
%\centerline {\includegraphics[width=5in]{Figure1.jpg}}

%\centerline{Figure 1}
%\vspace{.2in}
%
Suppose for a contradiction that $d(P,Q')\geq t'-t+J-(m+1)/q+1/q$. We will show that this implies the broken geodesics $PP_2Q'$ and $PQ_1Q'$ form a $K_0/2$-fat $1/q$-bigon. Since $X$ contains no such bigons we will conclude that in fact $d(P,Q')\geq t'-t+J-(m+1)/q+1/q$. 

To see that  $PP_2Q'$ is a $1/q$-path we use the fact that  $d(P,P_2)+d(P_2,Q')= (t_2-t)+d(P_2,Q')$ and that, by our induction hypothesis, $d(P_2,Q')< t'-t_2+J-m/q+1/q$:
\begin{eqnarray*}d(P,P_2)+d(P_2,Q')&=& (R+t_2-R-t)+d(P_2,Q')\\
&\leq& (t_2-t)+t'-t_2+J-m/q+1/q\\
&=& t'-t+J-m/q+1/q\leq d(P, Q')+1/q
\end{eqnarray*}
The last inequality comes from our supposition that $d(P,Q')\geq t'-t+J-(m+1)/q+1/q$.

A similar argument shows that the broken geodesic $PQ_1Q'$ is also a $1/q$-path if $d(P,Q')\geq t'-t+J-(m+1)/q+1/q$.

Now let $\xi$ be a point on the path $PQ_1Q'$ which minimises the distance to $P_2$. In particular $d(\xi, P_2)\leq d(Q_2, P_2)$. If $\xi$ lies on the arc $PQ_1$ then, by the triangle inequality, we compute:
\begin{eqnarray*}t_2-t&=&d(P,P_2)\leq d(P,\xi)+d(\xi,P_2)\leq d(P,Q_1)+d(Q_2,P_2)\\
&<&t_1-t+J-m/q+1/q+J+1/q.
\end{eqnarray*}
It follows that  $qJ=d_J(t_1,t_2)\leq t_2-t_1<2J+2/q$, which is a contradiction since we assumed $q\geq 3$ and $1\leq K_0\leq J$. 

It follows that  $\xi$ must lie on the arc $Q_1Q'$. But applying  Lemma \ref{synchronous} to the geodesics $\gamma, \gamma'$ with $s=R+t_2$ we see that $d(\xi, P_2)\geq K_0/2$ and hence the bigon is $K_0/2$-fat as required. Hence  $d(P,Q')< t'-t+J-(m+1)/q+1/q$ completing the induction.

Now by the triangle inequality
\begin{eqnarray*}R+t'&\leq&d(e,P)+d(P,Q')<R+t+(t'-t+J-n/q+1/q)\\
&=&R+t'+J-n/q+1/q\end{eqnarray*}
Hence $J-n/q+1/q\geq 0$, so $n\leq qJ+1$. So, $d_J(t,t')<(2^{qJ+1}-1)qJ$, and $r_0\leq (q(K_1-K_0)+1)(2^{qK_1+1}-1)qK_1+1$ because there are at least $r_0-1$ integer points in $[0, r_0)$ and each of these is a $J$-point for one of the  $q(K_1-K_0)+1$ possible values of $J$. 
\end{proof}
%%%%%%%%%%%%%%%%%%%%%%%%%%%%%%%%%%%%%%%%%%%%%%%%%%%%%%%%%%%%%%%%%%%%%%%%%%%%%%%%%%%%%%
%%%%%%%%%%%%%%%%%%%%%%%%%%%%%%%%%%%%%%%%%%%%%%%%%%%%%%%%%%%%%%%%%%%%%%%%%%%%%%%%%%%%%%
\section{Divergence functions and hyperbolicity}\label{div}
Recall that for a geodesic metric space $(X,d)$, a \emph{divergence function} is a map $f:\R_+\to\R$ such that for all $x\in X$, all $R\in\R_+$ and all geodesics $\gamma=[x,y]$, $\gamma'=[x,z]$ such that $d(\gamma(R),\gamma'(R))\geq f(0)>0$, if $r>0$ is such that $R+r\leq\min\{d(x,y),d(x,z)\}$ and $\alpha$ is a path in the closure of $X\setminus B(x,R+r)$ from $\gamma(R+r)$ to $\gamma'(R+r)$, then  the length of $\alpha$ is at least $f(r)$. We say that \emph{ geodesics diverge in $X$} if there is a divergence function $f$ so that $\lim_{r\to\infty}f(r)=\infty$. Papasoglu showed in \cite{Papa}[Corollary 1.3] that a geodesic metric space $(X,d)$ is hyperbolic if and only if geodesics diverge in $X$. Here we provide a quantitive version of this result in order to relate the hyperbolicity constant for a space to the eccentricity bound.

First, following \cite{Papa} we provide candidates for a divergence function.
Let $D>0$ and for $r>0$ define
$$f_D(r)=\inf\{d(\gamma(R+r),\gamma'(R+r))\hbox{ s. t. }\gamma(0)=\gamma'(0), d(\gamma(R),\gamma'(R))\geq D\}$$
Since $X$ is a geodesic space any path joining two points $\gamma(R+r), \gamma'(R+r)$ must have length at least $d(\gamma(R+r), \gamma'(R+r))$, and setting $f_D(0)=D$ it is easy to see that  the function $f_D$ is a divergence function for $X$.
\begin{propn}\label{thinBigonsquant}Let $(X,d)$ be a metric space such that any $q$-bigon is $4(q+\epsilon)$-thin. Then for $D>32/3+48\epsilon$ and $T>D/4-8\epsilon$, any $r_0>0$ such that $f_D(r_0)\leq T$ satisfies
$$r_0\leq (12T+26\epsilon-3D/4+1)(2^{12T+24\epsilon+1}-1)(12T+24\epsilon)+1.$$
In words, $\lim_{r\to\infty}f_D(r)=\infty$ and hence geodesic diverge.
\end{propn}
Before starting with the proof we show an intermediate result.
\begin{lemma} Suppose that $(X,d)$ is a metric space such that $q$-bigons are $K(q+\epsilon)$-thin for some constant $K\geq 1$ and some $\epsilon\geq 0$. Let $D\geq 2K\epsilon$ and $T\geq D/K-2\epsilon$. If two geodesics $\gamma$ and $\gamma'$ starting at the same point are such that $d(\gamma(R),\gamma'(R))\geq D$ for some $R\geq 0$ and $d(\gamma(R+r_0),\gamma'(R+r_0))=T$ for some $r_0>0$, then
$$d(\gamma(R+r),\gamma'(r+R))\in[D/K-2\epsilon,KT+2K\epsilon]$$
for any $r\in[0,r_0]$.
\end{lemma}
\begin{proof}Let us denote by $A_r:=d(\gamma(R+r),\gamma'(r+R))$. Take $a_r$ to be the midpoint on a geodesic $\alpha_r$ from $\gamma(R+r)$ to $\gamma'(R+r)$. Then the broken geodesics $\gamma(0)\gamma(R+r)a_r$ and $\gamma(0)\gamma'(R+r)a_r$ form the sides of an $A_r/2$-bigon, which by assumption is $K(A_r/2+\epsilon)$-thin. According to Lemma \ref{synchronous}, this bigon is at least $D/2$-fat since the distance between the two geodesics is more than $D$ at time $R$. Hence $A_r\geq D/K-2\epsilon$.

The upper bound is obtained in a similar way by looking at the broken geodesics $\gamma(0)\gamma(R+r)a_{r_0}$ and $\gamma(0)\gamma'(R+r)a_{r_0}$. They form a $T/2$-bigon, which by assumption is $K(T/2+\epsilon)$-thin. But by Lemma \ref{synchronous}, this bigon is at least $A_r/2$-fat since the distance between the two geodesics is $A_r$ at time $r\leq r_0$. Hence $A_r\leq KT+2K\epsilon$ as required.\end{proof}
We now can easily prove Proposition \ref{thinBigonsquant}.
\begin{proof}[Proof of Proposition \ref{thinBigonsquant}] Since $f_D(r_0)\leq T$  there are two geodesics $\gamma$ and $\gamma'$ and $R\in\R_+$ such that
\begin{itemize}
\item[(a)]$\gamma(0)=\gamma'(0)$
\item[(b)]$d(\gamma(R),\gamma'(R))\geq D$
\item[(c)]$d(\gamma(R+r_0),\gamma'(R+r_0))\leq T$.
\end{itemize}
We can apply Theorem \ref{meat} with $q=3$, $K_0=D/4-2\epsilon$ and $K_1=4T+8\epsilon$ which can be done using the previous lemma with $K=4$. The assumptions on $D$ and $T$ ensure that $K_1>K_0\geq 0$.\end{proof}
From now on we assume that $D>32/3+48\epsilon$ so that the function $f_D$ satisfies the conclusions of Proposition \ref{thinBigonsquant} for appropriate constants. The next step is to replace the divergence function $f_D$ by a divergence function $e$ of exponential growth. Given a rectifyable path $\alpha$, let us denote by $\ell(\alpha)$ its length. For $r>0$ we set 
\begin{eqnarray*}
e(r)=\inf_{R\in\R_+,\gamma,\gamma',d(\gamma(R),\gamma'(R))\geq D}\{\ell(\alpha)\mid \alpha \hbox{ a path from } \gamma(R+r) \hbox{ to }\\
 \gamma'(R+r) \hbox{ in } \overline{X\setminus B(x, R+r)}\hbox{ and } d(\gamma(R), \gamma'(R))\geq D\}\end{eqnarray*}
where $\gamma, \gamma'$ are geodesics, $\gamma(0)=\gamma'(0)=x$ and the infimum is taken over all geodesics $\gamma, \gamma'$ and all points $x\in X$ and all $R\in \mathbb R^+$.

It is clear that if we define $e(0)=f_D(0)=D$ then $e$ is a divergence function on $X$ and that $f(r)\leq e(r)$ for all $r\geq 0$. The following shows that this divergence function has exponential growth.
\begin{propn}\label{exponentialdivergence}For any $k>1$ and $r>u+kN$, then
$$e(r)>(3/2)^k(4N+2)$$
where $N=1+3D+\sup\{r|f_D(r)<9D\}$ and $u=\sup\{t|f_D(t)<4N+2\}$. In particular, $e$ has exponential growth.\end{propn}
We start with an intermediate result.
\begin{lemma}\label{smallLemma} If $\gamma$ and $\gamma'$ are geodesics with $\gamma(0)=\gamma(0)$ and $R>0$ satisfies $d(\gamma(R),\gamma'(R))\geq D$ then $d(\gamma(R+N),\gamma'(R+N))\geq 3D$, where $N$ is as in Proposition \ref{exponentialdivergence}.
\end{lemma}
\begin{proof}
Let $\beta$ be a geodesic joining $\gamma(R+N)$ to $\gamma'(R+N)$. Suppose for a contradiction that  $\beta$ has length less than $3D$ so no point on $\beta$ lies in the interior of the ball $B(x,R+N-3D)$, where $x=\gamma(0)$. Concatenate $\beta$ with the terminal subarcs of $\gamma, \gamma'$ of length $3D$ to obtain a path joining $\gamma(R+N-3D)$ to $\gamma'(R+N-3D)$. No point on this path lies in the interior of the ball $B(x,R+N-3D)$ so it must have length at least $e(N-3D)$. On the other hand we see that the path has length $6D+\ell(\beta)$ which by assumption is less than $9D$ so we get $e(N-3D)<9D$. However $e(N-3D)\geq f(N-3D)\geq 9D$ by choice of $N$ and this is a contradiction.
\end{proof}
\begin{proof}[Proof of Proposition \ref{exponentialdivergence}] Arguing by induction on $k$ it is enough to show that $e(r)\geq \frac{3}{2}e(r-N)$ for any $r>u+N$. 

%Indeed, we then by induction on $k$ see that for $r>u+kN$, $e(r)\geq {3\over 2}^ke(r-kN)$. But $r-kN>u$ so $e(r-kN)\geq 4N+2$, hence $e(r)\geq {3\over 2}^k(4N+2)$ as required.

%So it remains to show that $e(r)\geq {3\over 2}e(r-N)$ for any $r>u+N$.

 If $e(r)=\infty$ we are done so we can assume that $e(r)$ is finite and therefore there are geodesics $\gamma, \gamma'$  with $\gamma(0)=\gamma'(0)=x$ and $R>0$ such that $d(\gamma(R),\gamma'(R))\geq D$ such that there is an arc $\alpha$ in the closure of the complement of the ball $B(x,R+r)$ which joins $\gamma(R+r)$ and $\gamma'(R+r)$ and which has length less than $e(r)+1$.

Let $t_1=\sup\{t\in [0, M/2]\mid \alpha(t)\in B(x,R+r+N)\}$ and $t_2=\inf\{t\in [M/2, M]\mid \alpha(t)\in B(x, R+r+N)\}$.  Note that  $t_1\geq N$ and $M-t_2\geq N$ so the subarc $\alpha|_{[t_1, t_2]}$ must have length less than or equal to $M-2N$. Let $c_1$ and $c_2$ be geodesics from $x$ to $\alpha(t_1)$ and $\alpha(t_2)$ respectively. By combining the triangle inequality  lemma \ref{smallLemma} we see that 
\begin{eqnarray*}
3D&\leq &d(\gamma(R+N), \gamma'(R+N))\\
&\leq &d(\gamma(R+N), c_1(R+N))+d(c_1(R+N), c_2(R+N))+d(c_2(R+N), \gamma'(R+N)).
\end{eqnarray*}
It follows that at least one of the three terms in the above sum must be greater than or equal to $D$. It cannot be the middle term for the following reason: if $d(c_1(R+N),c_2(R+N))\geq D$ then, by definition of $e$, any path joining $c_1(R+N+r)$ to $c_2(R+N+r)$ which lies in the closure of the complement of the ball $B(x, R+N+r)$ must have length at least $e(r)$ so in particular the subarc $\alpha|_{[t_1, t_2]}$ must have length at least $e(r)$. On the other hand we observed that this subarc has length at most $M-2N$ so we see that $e(r)\leq M-2N$. Since $M\leq e(r)+1$ we have that $e(r)\leq e(r)+1-2N$ and since $N\geq 1$ this is a contradiction.

So we may assume (interchanging $\gamma$ and $c_1$ with $\gamma'$ and $c_2$ if necessary) that $d(\gamma(R+N), c_1(R+N))\geq D$. We construct a path $\alpha'$ in the complement of the ball $B(x, R+r)$ joining $\gamma(R+r)$ to $c_1(R+r)$ by concatenating the subarc $\alpha|_{[0, t_1]}$ with the subarc $c_1'$ of $c_1$ joining $c_1(R+r)$ to $\alpha(t_1)$. Since $\alpha(t_1)$ is within $R+r+N$ of $x$ we see that $c_1'$ has length at most $N$, and so $\alpha'$ has length at most $M/2+N$. On the other hand the length of $\alpha'$ is bounded below by $e(R+r-(R+N))=e(r-N)$, so we see that %
\[e(r-N)< M/2+N.\]
But recall that $M\leq e(r)+1$ so we see that $e(r)> 2e(r-N)-2N-1$. Finally, since  $r-N>u$ and by definition of $u$, we have that $e(r-N)\geq f_D(r-N)\geq 4N+2$. It follows that $e(r-N)/4\geq N+1/2$ and so 
\[e(r)> 2(e(r-N)-N-1/2)\geq 2e(r-N)-e(r-N)/4=\frac{3}{2}e(r-N).\]

\end{proof}

It follows from Lemma \ref{exponentialdivergence} that for any affine function $g(r)=ar+b$ there is some $r_0$ such that $e(r)>g(r)$ for all $r\geq r_0$. The value $r_0$ depends only on the function $f_D$ and the constants $a$ and $b$, though there may be no closed formula to compute it. It is clear that the value of $r_0$ may be bounded in terms of  the values of the constants $N$ and $u$ appearing in Lemma \ref{exponentialdivergence}, and those depend only on the constant $D$ and  the eccentricity bound $\epsilon$ as shown in the next Lemma. This will enable us to show that there is an upper bound on the hyperbolicity constant for $X$ which is a function of $\epsilon$ alone.
\begin{lemma}\label{constants}
Let $(X, d)$ be a geodesic metric space such that for each $q$ the $q$-bigons are $4q+4\epsilon$-thin. Let $D>32/3+48\epsilon$, $N=\sup\{r\mid f(r)<9D\}+1+3D$ and $u=\sup\{t\mid f(t)<4N+2\}$ as in Proposition \ref{exponentialdivergence}. Then $N$ and $r$ are bounded above by functions of $\epsilon$ and $D$. More precisely,
\begin{enumerate}
\item $N<(106D+26\epsilon+1)(2^{108D+24\epsilon+1}-1)(108D+24\epsilon)+2+3D$
\item $u\leq (48N+26\epsilon-3D/4+25)(2^{48N+25\epsilon+13}-1)(48N+24\epsilon+24)+1$
\end{enumerate}\end{lemma}
\begin{proof}
(1) Apply Proposition \ref{thinBigonsquant} with $T=9f_D(0)=9D$.

(2) Apply Proposition \ref{thinBigonsquant} with $T=4N+2>3D$.\end{proof}
\begin{thm}\label{thinBigons}Let $(X, d)$ be a geodesic metric space such that for each $q$ the $q$-bigons are $4q+4\epsilon$-thin. Then $X$ is $\delta(\epsilon)$-hyperbolic, for some function $\delta$ depending on $\epsilon$ alone.\end{thm}
\begin{proof}
%\vspace{.2in}
%\centerline {\includegraphics[width=15cm]{Figure2.jpg}}

%\centerline{Figure 2}
%\vspace{.2in}
%
To ensure that $D>32/3+48\epsilon$ we set $D=11+48\epsilon$. Let $x,y,z\in X$ and choose geodesics $\alpha_z=[x,y], \alpha_x=[y,z], \alpha_y=[z, x]$ in $X$. We denote by $\alpha_p^{-1}$ the reverse of the geodesic $\alpha_p$. We wish to estimate the thickness $\delta$ of this geodesic triangle $\Delta$ using the exponential divergence function $e$ defined above. The following argument is adapted from that given in \cite{short}.

Let $T_x=\sup\{t | d(\alpha_z(t),\alpha_y^{-1}(t))\leq D\}$ and set $x_y=\alpha_z(T_x)$ and  $x_z=\alpha_y^{-1}(T_x)$. Similarly we define $T_y$, $T_z$, $y_x$, $y_x$, $z_x$ and $z_y$. Now set $L_z=d(x,y)-(T_x+T_y)$, $L_x= d(y,z)-(T_y+T_z)$ and $L_y= d(z,x)-(T_z+T_x)$

CASE 1.  At least one of the values, say $L_z$, is non-positive. In this case any point on $\alpha_z$ is within $D$ of the other two sides, while any point on the subarcs $yy_z$ or $xx_z$ is within $D$ of the other two sides. There is $z'$ on $\alpha_z$ between $y_x$ and $x_y$ that is within $D$ of both $y_z$ and $x_z$. Hence the broken paths $z'y_zz$ and $z'x_zz$ form a $D$-bigon, consequently the triangle is  $4D+4\epsilon$-thin.

CASE 2. All three of the values $L_x, L_y, L_z$ are positive. In this case we will show that we can bound  all three of these values by some uniform value $L$ given in terms of $D$ and $\epsilon$, and hence, setting $D=11+48\epsilon$, $L$ depends on $\epsilon$ alone. Once this is done it is clear that any point on any side of the triangle is within $L/2+11+48\epsilon$ of some point on one of the other two sides, so it remains to find this bound $L$.

We can assume that $L_x\leq L_y\leq L_z$. First note that if $L_x\leq 2D$ we can run the argument from case 1 with $2D$ in place of $D$, so we may assume that $2D<L_x\leq L_y\leq L_z$.

We claim that the interior of $B(x,T_x+L_z/2)$ doesn't intersect the geodesic $\alpha_x$. Indeed, since $d=(x,y)=T_x+T_y+L_z$, the interior of the balls $B(x,T_x+L_z/2)$ and $B(y,T_y+L_z/2)$ have disjoint intersection. Similarly, the balls $B(x,T_x+L_z/2)$ and $B(z,T_z+L_y-L_z/2)$ have disjoint interiors. But since $d(y,z)\leq T_z+T_y+L_x$, the arc $\alpha_x$ is contained in the union of $B(z,T_z+L_y-L_z/2)$ and $B(y,T_y+L_z/2)$.

Now let $p=\alpha_z(T_x+L_z/2)$ and $p'=\alpha^{-1}_z(T_x+L_z/2)$. The arcs $py_x$, $y_xy_z$, $y_zz_y$, $z_yz_x$ and $z_xp'$ are in the complement of the ball $B(x,T_x+L_z/2)$ (since they are either in the ball $B(y,T_y+L_z/2)$ or $B(z,T_z+L_y-L_z/2)$). So concatenating those arcs we obtain a path in the closure of the complement of the ball $B(x, T_x+L_x/2)$ of  length $L_x+L_y+2D\leq 2D+4L_z/2$. Applying the divergence function to the geodesics $\alpha_z, \alpha_y^{-1}$ emanating from the point $x$  we see that $e(L_z/2)\leq 2D+4L_z/2$.

Now choose an integer $k\geq 0$ so that $L_z/2\in (u+kN, u+(k+1)N]$ where $u, N$ are the constants estimated in Lemma \ref{constants}. If $k\leq 0$ then $L_z<u+N$ and, since $u$, $N$ depend only on $\epsilon$, we are done. If $k\geq 1$ we can apply  Theorem \ref{exponentialdivergence} to show that $e(r)>(3/2)^k(4N+2)$. It follows that $k$ can also be bounded above in terms of $\epsilon$ and $D$. As noted before we may choose $D=11+48\epsilon$  to obtain a bound on $k$ in terms of $\epsilon$ alone. We denote this bound by $k(\epsilon)$. This gives a bound on $L_x/2$ since $L_x/2\leq u+(k(\epsilon)+1)N$. This in turn bounds the fatness of the triangle as less than or equal to $D+2(u+(k(\epsilon)+1)N=11+48\epsilon+2(u+(k(\epsilon)+1)N)$

We have shown that any geodesic triangle is either $10D+4\epsilon$-thin, which since $D=11+48\epsilon$ means that the triangle is  $110+484\epsilon$-thin, or it is $u+N$-thin, or it is $11+48\epsilon+2(u+(k(\epsilon)+1)N)$-thin. In the second case and third case the constants $u, k(\epsilon), N$ can all be written in terms of $\epsilon$ alone. If we take $\delta(\epsilon)=\max\{110+484\epsilon, u+N, 11+48\epsilon+2(u+(k(\epsilon)+1)N)\}$ then we see that the space $(X,d)$ is $\delta(\epsilon)$-hyperbolic as required.\end{proof}
%%%%%%%%%%%%%%%%%%%%%%%%%%%%%%%%%%%%%%%%%%%%%%%%%%%
\section{Proof of the quasi-balls characterization}\label{proofs}
%%%%%%%%%%%%%%%%%%%%%%%%%%%%%%%%%%%%%%%%%%%%%%%%%%%%
The proof of Theorem \ref{qb} is a sequence of simple observations, combined with Theorem \ref{thinBigons}.
Recall that if $A,B$ are subsets of a metric space $(X,d)$, then the Hausdorff distance between $A$ and $B$ is given by
$$d_H(A,B)=\inf\{r|A\subset N_r(B), B\subset N_r(A)\},$$
where for $r\geq 0$, $N_r(A)$ is the $r$-neighborhood of $A$. It is not clear how having eccentricity less than a constant $\delta$ and being at Hausdorff distance less than $\delta$ to a ball are related in general, but in case of intersection of balls in a metric space those notions are equivalent. Our first observation in this section holds for any geodesic metric space and gives the interior radius of the intersection of two balls.
\begin{lemma}\label{general}Let $(X,d)$ be a geodesic metric space.
\begin{itemize}\item[(1)] For any $x,y\in X$ with $d(x,y)=d$ and $s,t\geq 0$, if the balls $B(x,s)$ and $B(y,t)$ are neither disjoint nor nested, then
$$B(c,r)\subseteq B(x,s)\cap B(y,t),$$
where $r=\frac{s+t-d(x,y)}{2}$ and $c$ is any point on any geodesic between $x$ and $y$, at distance $\frac{s-t+d}{2}$ from $x$.
\item[(2)]If $s,t<d$ and $B(\xi,R)\subseteq B(x,s)\cap B(y,t)$, then $R\leq s+t-d$.
\end{itemize}
\end{lemma}
\begin{proof}(1) 
If the balls $B(x,s)$ and $B(y,t)$ are neither disjoint nor nested, then we have $r=\frac{s+t-d}{2}\geq 0$ and $0<\frac{s-t+d}{2}<d$, and hence given any geodesic $\gamma$ from $x$ to $y$ we may take a point $c$ on $\gamma$ at distance $\frac{s-t+d}{2}$ from $x$ and a ball of radius $r$ around $c$. Then for $z\in B(c,r)$
$$d(x,z)\leq d(x,c)+d(c,z)\leq \frac{s-t+d}{2}+ \frac{s+t-d}{2}=s$$
and similarly $d(y,z)\leq t$.

(2) First notice that our assumptions on $s$ and $t$ being strictly smaller than $d$ show that $x$ and $y$ do not belong to $B(x,s)\cap B(y,t)$. Let $a=d(x,\xi)$ and $b=d(y,\xi)$, so that $a\leq s$, $b\leq t$ and $d\leq a+b\leq s+t$. Take a geodesic $\gamma_{x\xi}$ from $x$ to $\xi$ and a point $z$ on this geodesic, at distance $R$ from $\xi$. Such a point exists because $x$ does not belong to $B(\xi,R)$. Since $z\in B(\xi,R)\subseteq B(y,t)$ we have
$$d=d(x,y)\leq d(x,z)+d(z,y)\leq a-R+t.$$
Similarly, take a geodesic $\gamma_{y\xi}$ from $y$ to $\xi$ and a point $z'$ on this geodesic, at distance $R$ from $\xi$. Since $z'\in B(\xi,R)\subseteq B(x,s)$ we have
$$d=d(x,y)\leq d(x,z')+d(z',y)\leq b-R+s.$$
Combining the 2 inequalities gives $2d\leq a+b-2R+s+t$, hence $2R\leq s+t-d+a+b-d\leq 2(s+t-d)$.\end{proof}
We can now prove one implication of Theorem \ref{qb} and Corollary \ref{RtreesChar}.
\begin{propn}\label{hyp=>qb} If a geodesic metric space $(X,d)$ is $\delta$-hyperbolic with hyperbolicity constant less than or equal to $\delta\geq 0$, then both the eccentricity of the intersection of any two balls and the Hausdorff distance from the intersection to a ball are both uniformly bounded by $2\delta$.\end{propn}
\begin{proof}Take $x,y\in X$ with $d=d(x,y)$ and $s,t\in\R_+$, with $s\geq t$. We will show that the eccentricity of $B(x,s)\cap B(y,t)$ is less than $2\delta$. This implies the statement about Hausdorff distance as well, by definition
of Hausdorff distance (recalled in the beginning of this section).

According to Lemma \ref{general} part (1) either $B(x,s)\subseteq B(y,t)$, $B(y,t)\subseteq B(x,s)$, $B(x,s)\cap B(y,t)=\emptyset$ or $B(c,r)\subseteq B(x,s)\cap B(y,t)$, where $c$ and $r$ are as defined in that lemma. In the first three cases the eccentricity is clearly bounded by $0$. In the remaining case it suffices to show that there is a constant $\epsilon$ independent of $x,y,s,t$  such that $B(x,s)\cap B(y,t)$ is  contained in some ball of radius $r+\epsilon$. We will  show that in fact $B(c,r)\subseteq B(x,s)\cap B(y,t)\subseteq B(c,r+2\delta)$. Now, for $z\in B(x,s)\cap B(y,t)$, let us estimate the distance to $c$. 

Since $c$ lies on a geodesic from $x$ to $y$ it is within $\delta$ of a point $p$ which lies on a geodesic from $y$ to $z$ or on a geodesic from $x$ to $z$. We first assume that $p$ lies on a geodesic from $x$ to $z$. By the triangle inequality, we have that $d(x,c)+d(c,z)\leq d(x,p)+d(p,z)+2\delta$. Since $p$ lies on a geodesic from $x$ to $z$ this yields $d(x,c)+d(c,z)\leq d(x,z)+2\delta=s+2\delta$. Now we have, as required:
$$d(c,z)\leq s+2\delta-\frac{s-t+d(x,y)}{2}=\frac{s+t-d(x,y)}{2}+2\delta.$$
If $p$ lies on a geodesic from $y$ to $z$ instead, then we use the same argument switching the roles of $x,y$ and of $s,t$.\end{proof}
The following says that in a geodesic metric space such that the intersection of any two balls has uniformly bounded eccentricity, then the set of points on $K$-paths is uniformly close to a geodesic.
\begin{lemma}\label{dpathTravelClose}Suppose that $(X,d)$ has a uniform bound $\epsilon\geq 0$ on the eccentricity of the intersection of any two balls. Then for any $q\geq 0$, given any two points $x,y\in X$, any point on a $q$-path from $x$ to $y$ is contained in the $2q+2\epsilon$-neighbourhood of any geodesic from $x$ to $y$.
\end{lemma}
\begin{proof}Let $z$ be a point on a $q$-path $\mu$ from $x$ to $y$ with $s=d(x,z)$, $t=d(y,z)$ so that $d\leq s+t\leq d+q$. If $s\geq d$ then $t\leq q$ and $z$ is within $q$ of $y$. Similarly if $t\geq d$ then  $s\leq q$ and $t$ is within $q$ of $x$ so we may assume that both $s,t< d$. Now let $Y=B(x,s)\cap B(y,t)$, so that $z\in Y$ so by Lemma \ref{general} point (2) we see that any ball contained in $Y$ has radius at most $q$. It follows from the bounded eccentricity hypothesis that $Y\subseteq B(\xi, q+\epsilon)$ for some point $\xi\in X$. Let $\gamma$ be a geodesic from $x$ to $y$ and $c$ any point on $\gamma$ at distance less than $s$ to $x$ and less than $t$ to $y$ (such a point exists because $d\leq s+t$). Then $c,z\in Y\subseteq B(\xi,q+\epsilon)$, hence $d(c,z)\leq d(c,\xi)+d(\xi,z)\leq 2q+2\epsilon$.\end{proof}
\noindent
{\sl Remark: Taking $K=0$ in the lemma above shows that in a geodesic metric space such that the intersection of any two metric balls has eccentricity less than or equal to $\delta$, any geodesic between two points is contained in a $\delta$-neighbourhood of any other geodesic between those two points.}

An analogous result holds in terms of Hausdorff distance.
\begin{lemma}\label{dpathTravelClose2}Suppose that $(X,d)$ has a uniform bound $\epsilon\geq 0$ on the Hausdorff distance from the intersection of any two balls to a ball. Then for any $q\geq 0$, given any two points $x,y\in X$, any point on a $q$-path from $x$ to $y$ is contained in the $2q+6\epsilon$-neighbourhood of any geodesic from $x$ to $y$.
\end{lemma}
\begin{proof}Let $z$ be a point on a $q$-path $\mu$ from $x$ to $y$ with $s=d(x,z)$, $t=d(y,z)$ so that $d\leq s+t\leq d+q$. If $s+\epsilon\geq d$ then $t-\epsilon\leq q$ and $z$ is within $q+\epsilon$ of $y$. Similarly if $t+\epsilon\geq d$ then  $s-\epsilon\leq q$ and $t$ is within $q+\epsilon$ of $x$ so we may assume that both $s+\epsilon,t+\epsilon< d$. Now let $Y=B(x,s)\cap B(y,t)$, so that $z\in Y$. By assumption, $d(Y,B(\xi,R))\leq\epsilon$, for some $\xi\in X$ and some $R\geq 0$, which implies that $Y\subset B(\xi,R+\epsilon)=N_{\epsilon}(B(\xi,R))$ and that $B(\xi,R)\subset N_{\epsilon}(Y)\subset B(x,s+\epsilon)\cap B(y,t+\epsilon)$. So by Lemma \ref{general} point (2) we see that any ball contained in $B(x,s+\epsilon)\cap B(y,t+\epsilon)$ has radius at most $q+2\epsilon$. It follows that $Y\subseteq B(\xi, q+3\epsilon)$. Let $\gamma$ be a geodesic from $x$ to $y$ and $c$ any point on $\gamma$ at distance less than $s$ to $x$ and less than $t$ to $y$ (such a point exists because $d\leq s+t$). Then $c,z\in Y\subseteq B(\xi,q+3\epsilon)$, hence $d(c,z)\leq d(c,\xi)+d(\xi,z)\leq 2q+6\epsilon$.\end{proof}
We can now prove the other implication in Theorem \ref{qb}, namely that (b) implies (a). 
\begin{lemma}\label{qb=>hyp}There is a function $\delta:\mathbb R^+\rightarrow \mathbb R^+$ such that:
\begin{enumerate}
\item If $(X,d)$ is a geodesic metric space with the property that the intersection of any two balls has eccentricity bounded by $\epsilon$ then $X$ is $\delta(\epsilon)$-hyperbolic. 
\item If $(X,d)$ is a geodesic metric space with the property that the intersection of any two balls is at Hausdorff distance less than $\epsilon/3$ then $X$ is $\delta(\epsilon)$-hyperbolic.\end{enumerate}
\end{lemma}
\begin{proof}According to Lemma \ref{dpathTravelClose} or Lemma \ref{dpathTravelClose2}, in such a metric space all $(1,q)$ bigons are $4(q+\epsilon)$-slim. We conclude using Theorem \ref{thinBigons} above.\end{proof}
As Pierre Pansu pointed out to us it immediately follows that we can take $\delta$ to be linear in $\epsilon$.
\begin{cor}
For any $\epsilon>0$ the space $(X,d)$ is $\epsilon\delta(1)$ hyperbolic.
\end{cor}
\begin{proof}
Scaling the metric we see that $(X, d/\epsilon)$ has the property that the intersection of any two balls has eccentricity bounded by $1$ and so is $\delta(1)$-hyperbolic. Rescaling we see that $(X, d)$ is $\epsilon\delta(1)$-hyperbolic.
\end{proof}
We conclude this section with the proof of Corollary 2 which asserts that $\R$-trees are characterised by the property that  the intersection of any two metric balls is a metric ball.

\begin{proof}[Proof of Corollary \ref{RtreesChar}]One implication is given by Proposition \ref{hyp=>qb}. Conversely,  if the space $(X,d)$ has the property that the intersection of any two balls has eccentricity 0 then for any $\epsilon>0$ the intersection of any two balls has eccentricity bounded by $\epsilon$ and so the space is $\epsilon\delta(1)$ hyperbolic for all $\epsilon>0$. It is therefore $0$-hyperbolic and hence must be an $\mathbb R$-tree.\end{proof}
{\sl Remark: Notice that in fact we do not need to assume that the intersection of any two balls is a ball to carry out the proof, only that the eccentricity of such an intersection is $0$. A priori this is a weaker condition, however  in an $\mathbb R$-tree the intersection of two balls is always a ball and therefore, as a consequence of the theorem, the two conditions are equivalent.}
%%%%%%%%%%%%%%%%%%%%%%%%%%%%%%%%%%%%%%%%%%%%%%%%%%%%%%%%%%%%
%%%%%%%%%%%%%%%%%%%%%%%%%%%%%%%%%%%%%%%%%%%%%%%%%%%%%%%%
\section{Miscellaneous comments}\label{ext}
In Theorem \ref{qb} and Corollary \ref{RtreesChar}, the assumption that the metric space $(X,d)$ be geodesic might not be needed. The notion of hyperbolic spaces extends to non-geodesic metric spaces via the Gromov product (see e.g. D\'efinition 3, page 27 of \cite{GH}) and it would be interesting to find an appropriate generalisation of these results to that context. In particular recall that a $\delta$-ultrametric space is a metric space $(X,d)$ which satisfies the following strengthened version of the triangle inequality,
$$d(x,y)\leq\max\{d(x,z),d(z,y)\}+\delta$$
for all $x,y,z$ in $X$. It is easy to see that at least two of $d(x,y)$, $d(y,z)$ and $d(x,z)$ differ by at most $\delta$, meaning that any triangle is almost isoceles. These are examples of $2\delta$-hyperbolic spaces in the sense of Gromov, see \cite{Gr} Section 1.2  on page 90. It would be interesting to know if those spaces do satisfy the property that any intersection of two balls is almost a ball (in some sense).
\begin{defi} We say that a geodesic metric space $(X, d)$ has the geodesic extension property if any geodesic arc $\gamma:[0, a]\rightarrow X$ extends to a geodesic $\gamma':[0, \infty)\rightarrow X$, i.e., $\gamma'|_{[0, a]}=\gamma$.
\end{defi}
The most important feature of a space with the geodesic extension property is that centres and radii of balls are well defined (this is easily checked). More precisely we will use the following.
\begin{lemma}\label{outeradius}Let $(X,d)$ be a space with the geodesic extension property. For any $x,y\in X$ and $s,t\geq 0$, with $s,t\leq d(x,y)$. Then $B(c,r)$ is the biggest ball that fits in $B(x,s)\cap B(y,t)$, where $r=\frac{s+t-d(x,y)}{2}$ and $c$ is a point on any geodesic between $x$ and $y$, at distance $\frac{s-t+d}{2}$ from $x$.\end{lemma}
\begin{proof}Let $Y=B(x,s)\cap B(y,t)$. If $s+t< d(x,y)$ then $Y=\emptyset$ and there is nothing to prove. So let us assume that $s+t\geq d(x,y)$. Take $r\geq 0$ and $c\in Y$ such that $B(c,r)\subseteq Y$. Let $a=d(x,c)$ and $b=d(y,c)$, so that $a+b\geq d(x,y)$. Since $B(c,r)\subseteq B(x,s)$, we deduce that $a+r\leq s$, and similarly, since $B(c,r)\subseteq B(y,t)$, we deduce that $b+r\leq t$ (here we use the geodesic extension property). Combining those two inequalities shows that
$$2r\leq s+t-(a+b)\leq s+t-d(x,y).$$
\end{proof}
If $X=X_\kappa$ is the symmetric space of constant curvature $\kappa\leq 0$ there is a single triangle $x$, $y$, $z$ (up to isometry) with side lengths $s=d(x,z)$, $t=d(y,z)$ and $d=d(x,y)$, so we define 
$$\hbox{Ecc}_\kappa(s,t,d):=d(z,c)-\frac{s+t-d}{2},$$
where $c$ is the point on the geodesic between $x$ and $y$ at distance $(d+s-t)/2$ from $x$. (This point exists because the triangle inequality ensures that $t\leq d+s$ so the distance is positive, and because $s\leq d+t$ implies that $(d+s-t)/2\leq d$.)  We shall see that there is an analogue to Theorem \ref{qb} which characterises CAT($\kappa$) geometry for $\kappa\leq 0$
\begin{thm}\label{CAT0} Given $\kappa\leq 0$, a geodesic metric space $(X,d)$ with geodesic extension property is $CAT(\kappa)$ if and only if the eccentricity of the intersection of any two balls of respective radii $s$ and $t$ and at distance $d$ is bounded by $\hbox{Ecc}_\kappa(s,t,d)$.
\end{thm}
\begin{proof}One implication is clear so we suppose by contradiction that  the eccentricity of the intersection of any two balls of respective radii $s$ and $t$ and at distance $d$ is bounded by $\hbox{Ecc}_\kappa(s,t,d)$ but that $X$ is not CAT($\kappa$). Then there is a geodesic triangle $x,z,y$ in $X$ and a point $p=\gamma(r)$ on a geodesic $\gamma$ from $x$ to $y$ such that  $d(\overline x, \overline p)< d(x,p)$. Let $r_0=\sup\{r'<r\mid d(\gamma(r'), z)\leq d(\overline\gamma(r'), \overline z)\}$ where $\overline\gamma$ denotes the geodesic $\overline x\overline y$ in the comparison triangle. Similarly let $r_1=\inf\{r'>r\mid d(\gamma(r'), z)\leq d(\overline\gamma(r'), \overline z)\}$. In other words, $r_0$ and $r_1$ are the nearest points left and right of $p$ that satisfy the CAT($\kappa$) inequality. So for any $r'$ in the open interval  $(r_0, r_1)$ we have $d((\gamma(r'), z)>d(\overline\gamma(r'), \overline z)$. Notice that $r_1>r_0$ since the metric varies continuously with points. Set $s'=d(\gamma(r_0),z)=d(\overline\gamma(r_0),\overline z)$, and  
$t'=d(\gamma(r_1), z)=d(\overline\gamma(r_1), \overline z)$. Let $d'=d(\gamma(r_0), \gamma(r_1))=r_1-r_0$.

Now consider a geodesic triangle $\gamma(r_0),z,\gamma(r_1)$ where the geodesic from $\gamma(r_0)$ to $\gamma(r_1)$ is taken to be the restriction of $\gamma$ to the closed interval $[r_0, r_1]$. Clearly the geodesic triangle $\overline\gamma(r_0),\overline z,\overline\gamma(r_1)$ is a comparison triangle in $X_\kappa$.

%\vspace{.2in}
%   \includegraphics[width=4in]{Figure3.jpg}   
%   
%\centerline{Figure 3}
%\vspace{.2in}

Take the point $c'$ on the above geodesic from $\gamma(r_0)$ to $\gamma(r_1)$ at distance $(s'+d'-t')/2$ from $\gamma(r_0)$. (This point exists because of the triangle inequality.) First notice that $c'$ has to be equal to either $\gamma(r_0)$ or $\gamma(r_1)$, otherwise we would have $d(z, c')>d(\overline z, \overline c')$ and  the intersection of $B(\gamma(r_0),s')$ with $B(\gamma(r_1),t')$ exceeds the allowed eccentricity. Assume that $c'=\gamma(r_0)$ (if $c'=\gamma(s_1)$ the argument is similar and we omit it). Then $(s'+d'-t')/2=0$ so $t'=s'+d'$. It follows that the path given by concatenating the geodesics from $z$ to $\gamma(r_0)$ and from $\gamma(r_0)$ to $\gamma(r_1)$ is itself a geodesic. Now let $m$ be the midpoint of the geodesic from $\gamma(r_0)$ to $\gamma(r_1)$. This point is at distance $s'+(r_1-r_0)/2$ from $z$. Inspecting the (degenerate) comparison triangle we see that $\overline m$ is also at distance $s'+(r_1-r_0)/2$ from $\overline z$ but this contradicts our assumption that every point between $\gamma(r_0)$ and $\gamma(r_1)$ is further from $z$ than we see in the comparison triangle. 
\end{proof}

\sl{Remark:
Following the arguments in the proof of Theorem 1 we note that there is an alternative notion of a "based" eccentricity function which measures eccentricity from the defined centre $c$ for both the inscribed and circumscribed balls. In these terms it is easy to see that Gromov hyperbolicity is equivalent to the existence of a uniform bound on the based eccentricity function, while the proof of Theorem \ref{CAT0} shows that CAT($\kappa$) geometry is characterised by bounding the based eccentricity function in terms of the function ${Ecc}_\kappa$. Hence both notions of non-positive curvature may be naturally expressed in terms of eccentricity bounds.}

\vfill\eject

%%%%%%%%%%%%%%%%%%%%%%%%%%%%%%%%%%%%%%%%%%%%%
\section{Appendix: Pomroy's result}
%%%%%%%%%%%%%%%%%%%%%%%%%%%%%%%%%%%%%%%%%%%%%
In \cite{Papa} Papasoglu showed that for a graph hyperbolicity was equivalent to a bound on the thinness of geodesic bigons. As remarked before the same statement is not true for general geodesic metric spaces (any non-hyperbolic CAT(0) space furnishes a counter example since uniqueness fo geodesics gives a bound of $0$ on the fatness of geodesic bigons.) The point is that the bound on the fatness of geodesic bigons in a graph gives an automatic bound on the fatness of $(1,1)$ quasi-geodesic bigons and Papasoglu remarks that there is a natural generalisation of the result as follows. The theorem appears in the Masters dissertation of Pomroy \cite{Pomroy} but to the best of our knowledge no proof exists in the literature. We offer a proof of the result in order to place it on the record.
\begin{thm}[Pomroy \cite{Pomroy}]Let $(X,d)$ be a geodesic metric space. If there is $\epsilon,\rho>0$ so that $\rho$-bigons are uniformly $\epsilon$-thin, then $X$ is hyperbolic.\end{thm}
\begin{proof}We argue by contradiction, similarly to Corollary \ref{thinBigons} and using those particular divergence functions. According to Theorem \ref{thinBigonsquant} for any $D\geq 0$, the functions $f_D$ as defined above never tend to infinity, and hence there is $L=L(D)$ so that $\lim\inf f_D=L(D)<\infty$. This means that for any $t_0\in\R_+$, there are two geodesics $\gamma$ and $\gamma'$ and $R\in\R_+$ such that
\begin{itemize}
\item[(a)]$\gamma(0)=\gamma'(0)$
\item[(b)]$d(\gamma(R),\gamma'(R))\geq D$
\item[(c)]$d(\gamma(R+r_0),\gamma'(R+r_0))\leq L+1$ for some $r_0\geq t_0$.
\end{itemize}
Let us fix $q\geq 3$ so that $1/q\leq\rho$ and again let $A_r=d(\gamma(R+r),\gamma'(R+r))$ and $a_0=\gamma(R+r_0),a_1,\dots,a_n=\gamma'(R+r_0)$ be points on a geodesic from $\gamma(R+r_0)$ to $\gamma'(R+r_0)$ and at distance less than $\rho$ from each other (we can choose $n\leq (L+1)q+1$). For $1=1,\dots,n$, we see a $\rho$-bigon as follows: one side is a geodesic from $\gamma(0)$ to $a_i$, and the other is a broken geodesic $\gamma(0)a_{i-1}a_i$. So our assumptions say that it is $\epsilon$-thin, and hence $A_r\leq n\epsilon$. Now let $b_0=\gamma(R+r),b_1,\dots,b_m=\gamma'(R+r)$ points on a geodesic from $\gamma(R+r)$ to $\gamma'(R+r)$ and at distance less than $\rho$ from each other (we can choose $m\leq A_rq+1$). Again we construct a $\rho$-bigon as follows: one side is a geodesic from $\gamma(0)$ to $b_i$, and the other is a broken geodesic $\gamma(0)b_{i-1}b_i$. Moreover, since $A_0=D$, it means that $m\epsilon\geq D$, and hence $(A_rq+1)\epsilon\geq D$, so that $A_r\geq (D/\epsilon-1)1/q$. This means that the geodesics $\gamma$ and $\gamma'$ fulfill the assumptions of Theorem \ref{meat} for any $D$ big enough (i.e., $D$ so that $(D/\epsilon-1)1/q>0$, and hence,  taking $D$ big enough (i.e., so that $(D/\epsilon-1)1/2q>\epsilon$), Theorem \ref{meat} contradicts our assumption that $\rho$-bigons are uniformly $\epsilon$-thin.
\end{proof}
\noindent
{\sl Remark: This may seem close to Lemma 7.2 in \cite{MM} which states that if $(3,0)$-quasigeodesics stay uniformly close to any geodesic between the endpoints, then the space is hyperbolic. However  the proof of Papasoglu's or Pomroy's result is considerably more elaborate.}
\vfill\eject

   \includegraphics[width=5in]{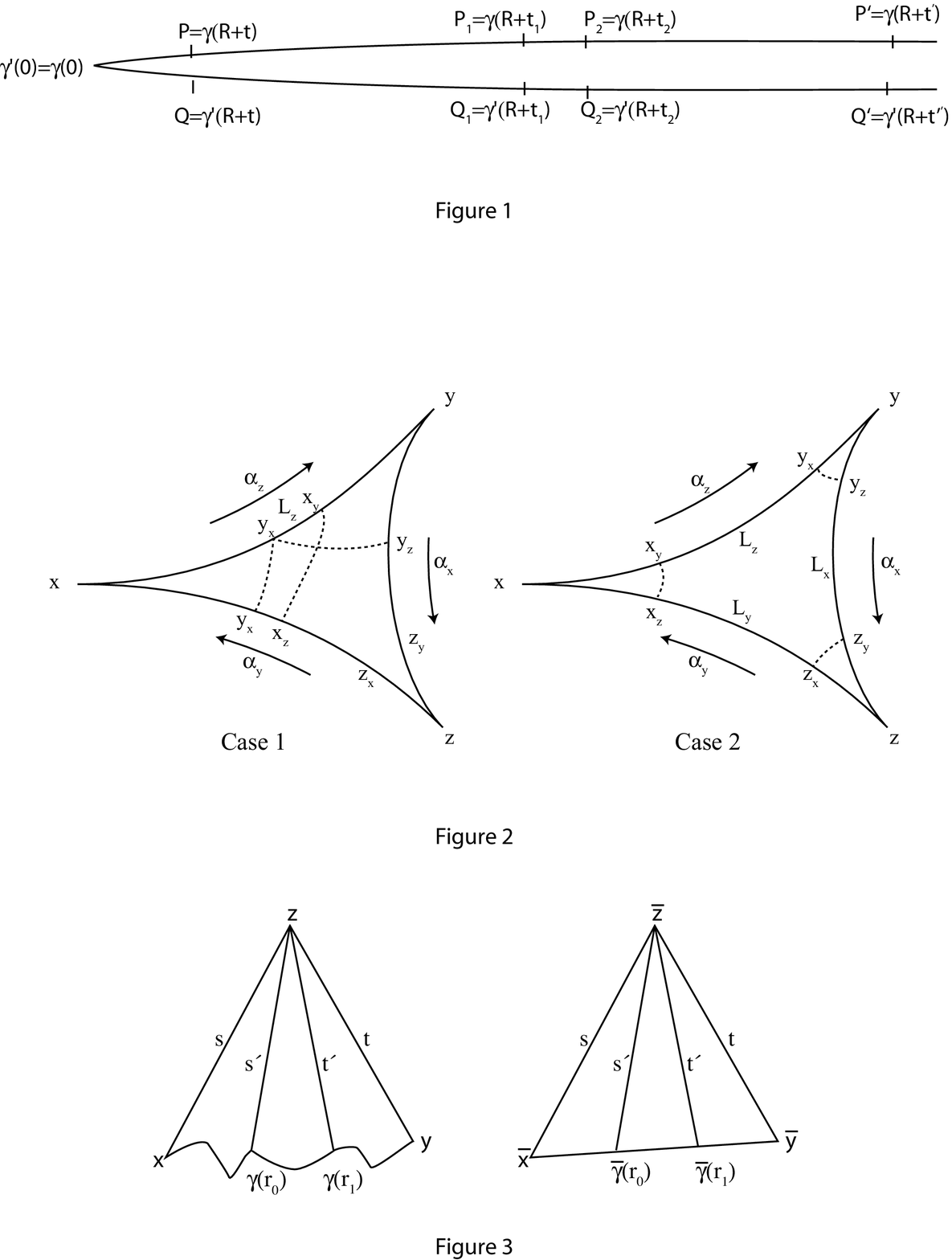}   
   \eject

 %%%%%%%%%%%%%%%%%%%%%%%%%%%%%%%%%%%%%
\end{document}